\documentclass[11pt]{amsart}
\usepackage{geometry}                % See geometry.pdf to learn the layout options. There are lots.
\usepackage{graphicx}
\usepackage{amssymb}
\usepackage{amscd}
\usepackage{epstopdf}
\newtheorem{df}{Definition}[section]

\newtheorem{thm}{Theorem}[section]
\newtheorem{prop}{Proposition}[section]
\newtheorem{lm}{Lemma}[section]
\newtheorem{ex}{Example}[section]
\newtheorem{remark}{Remark}[section]
\newtheorem{fact}{Fact}[section]
\newtheorem{cor}{Corollary}[section]

\title{On a generalization of Deuring's results}
\author{Ken-ichi Sugiyama}

\begin{document}
\maketitle

\begin{center}
Department of Mathematics and Informatics, Faculty of Science,\\ 
Chiba University, 1-33 Yayoi-cho Inage-ku,\\
Chiba 263-8522, Japan \\
e-mail address : sugiyama@math.s.chiba-u.ac.jp
\end{center}

\begin{abstract} Using the Dieudonn\'{e} theory we will study a reduction of an abelian variety with complex multiplication at a prime. Our results may be regarded as generalization of the classical theorem due to Deuring for CM-elliptic curves. We will also discuss a sufficient condition for a prime at which the reduction of a CM-curve is maximal.
AMS classification 2010: 11G10, 11G15, 11G20, 11G30, 14F30, 14F40, 14G50, 14H40.
\end{abstract}
\section{Introduction}
Let $E$ be an elliptic curve over the field of rational numbers ${\mathbb Q}$ with complex multiplication (which will be abbreviated by {\it CM}) by the integer ring of an imaginary quadratic field $K$ that has a non-singular reduction $E_{p}$ at a prime $p$. The classical Deuring's theorem states that it is ordinary or supersingular according to whether $p$ completely splits or remains prime in $K$, respectively. It is known that $E_{p}$ is ordinary (resp. supersingular) if and only if the $p$-divisible group $E_{p}[p^{\infty}]$ is ${\mathbb L}$ (resp. $G_{1.1}$) (see \$2 for details). Suppose that $E_{p}$ is supersingular. If $p$ is greater than $3$ the characteristic polynomial of the $p$-th power Frobenius on an $l$-adic Tate module ($l \neq p$) is $t^{2}+p$ and the number of ${\mathbb F}_{p^2}$-rational points attains the Hasse-Weil upper bound:
\[|E_{p}({\mathbb F}_{p^2})|=1+p^2+2p,\]
where $|\cdot|$ stands for cardinality. In this report we will generalize these results to a proper smooth curve over ${\mathbb Q}$ of a higher genus with CM. \\

Let us first explain terminologies used in the paper. Let $X$ be an object defined over a field $F$. The base change over an extension $F^{\prime}$ of $F$ is denoted by $X\otimes_{F}F^{\prime}$. If $F^{\prime}=\overline{F}$, the algebraic closure, it is simply described by $\overline{X}$. Let $V$ be a proper smooth variety defined over ${\mathbb F}_q$ ($q=p^{f}$). Then 
 $\Phi_{V,q}(t)$ stands for the characteristic polynomial of the $q$-th power Frobenius on $H^{1}_{et}(\overline{V},{\mathbb Q}_{l})$ ($l\neq p$).\\

%An algebraic closure of a field $F$ is written by $\overline{F}$. For an object $X$ defined over $F$ its base change to an field extension $F^{\prime}$ of $F$ is denoted by $X\otimes_{F}F^{\prime}$. For simplicity we sometimes write $\overline{X}:=X\otimes_{F}\overline{F}$. 
% 
% Suppose that $F$ is a finite field of $q$-elements and let $V$ be a proper smooth variety. Then the characteristic polynomial of the $q$-power Frobenius on the first $l$-adic \'{e}tale cohomology group $H^{1}_{et}(\overline{V},{\mathbb Q}_{l})$ ($l\neq {\rm char}(F)$) is described by $\Phi_{V,q}(t)$.\\

Let $k$ be a field of characteristic $p$.  An abelian variety $A$ of dimension $g$ over $k$ will be called {\it supersingular} (resp. {\it superspecial}) if $\overline{A}$ is isogeneous (resp. isomorphic) to a product of supersingular elliptic curves and {\it ordinary} if the group of $p$-torsion points $\overline{A}[p]$ is isomorphic to $({\mathbb Z}/p{\mathbb Z})^{g}$. We mention that $A$ has {\it CM} (or sometimes {\it CM by} ${\mathcal O}_{K}$) if there is a finite extension $k^{\prime}$ of $k$ such that  the endomorphism ring of $A\otimes_{k}k^{\prime}$ contains the integer ring ${\mathcal O}_{K}$ of a CM-field $K$ satisfying $[K:{\mathbb Q}]=2g$. %Let $W(k)$ be the ring of Witt vectors whose coefficients in $k$ and ${\mathcal A}$ an abelian scheme over $W(k)$ of relative dimension $g$. We call ${\mathcal A}$ has {\it CM} (or sometimes {\it CM by} ${\mathcal O}_{K}$) if the endomorphism ring of ${\mathcal A}\otimes_{W(k)}W(k^{\prime})$ contains the integer ring ${\mathcal O}_{K}$ of a CM-field $K$ satisfying $[K:{\mathbb Q}]=2g$, where $k^{\prime}$ is a finite extension of $k$. Let $A$ be an abelian variety over $k$ with CM by ${\mathcal O}_{K}$. We say it {\it lifts} to $W(k)$ if there is an abelian scheme ${\mathcal A}$ over $W(k)$ with CM by ${\mathcal O}_{K}$ such that $(A,{\mathcal O}_{K})$ is the reduction of $({\mathcal A}, {\mathcal O}_{K})$.\\
Let $K_{0}$ be the maximal totally real subfield and $[K_{0}:{\mathbb Q}]=g$. We assume that $p$ is unramified in $K$ and let
\begin{equation}
p={\frak P}_{1}\cdots {\frak P}_{t}
\end{equation}
be the prime factorization in $K_{0}$.

%Let  $W(k)$ be the ring of Witt vectors whose coefficients in $k$. An abelian scheme ${\mathcal A}$ over $W(k)$ of relative dimension $g$ will be called {\it with CM} (or sometimes {\it with CM by} ${\mathcal O}_{K}$) if the endomorphism ring over $W(\overline{k})$ contains the integer ring ${\mathcal O}_{K}$ of a CM-field $K$ satisfying $[K:{\mathbb Q}]=2g$.  Let $K_{0}$ be the totally real subfield with $[K_{0}:{\mathbb Q}]=g$. We assume that $p$ is unramified in $K$ and let
%\begin{equation}
%p={\frak P}_{1}\cdots {\frak P}_{t}
%\end{equation}
%be the prime factorization in $K_{0}$.
\begin{thm}
Let $A$ be an abelian variety of dimension $g$ defined over a finite field $k$ of characteristic $p$ endowed with CM by ${\mathcal O}_{K}$.
%Let $k$ be a finite field of characteristic $p$ and $A$ an abelian variety of dimension $g$ defined over $k$ with CM by ${\mathcal O}_{K}$. 
Suppose that every ${\frak B}_{i}$ of (1) remains prime in $K$. Then $A$ is supersingular. If moreover $t=g$ (i.e. $p$ completely splits in $K_{0}$) it is superspecial.
\end{thm}
If $p \geq 5$ this determines the characteristic polynomial of $A$. In fact suppose that the assumptions of {\bf Theorem 1.1} are satisfied. Then there is a positive integer $m$ so that $A\otimes_{k}{\mathbb F}_{p^m}$ is isogeneous to a product of supersingular elliptic curves $\{E_{i}\}_{1\leq i \leq g}$ defined over ${\mathbb F}_{p^m}$ and $\Phi_{A,p^m}(t)=\prod_{i=1}^{g}\Phi_{E_{i},p^m}(t)$. Since $p\geq 5$, $\Phi_{E_{i},p^m}(t)$ is a one of the following (\cite{Waterhouse1969} {\bf Theorem 4.1}):
\[t^{2}+p^{m},\quad t^{2}\pm p^{\frac{m}{2}}t+p^{m}, \quad t^{2}\pm 2p^{\frac{m}{2}}t+p^{m},\]
where the last two occur when $m$ is even.
\begin{thm}
Let $A$ be an abelian variety of dimension $g$ defined over a finite field $k$ of characteristic $p$ endowed with CM by ${\mathcal O}_{K}$.
%Let $k$ be a finite field of characteristic $p$ and $A$ an abelian variety of dimension $g$ defined over $k$ with CM by ${\mathcal O}_{K}$. 
Suppose that $p$ completely splits in $K$ then $A$ is ordinary.
\end{thm}
%\begin{remark} In these theorems we consider the assumption of a lift may be superfluous. It is necessary to show that the $p$-adic Tate module $T_{p}(A)$ is free ${\mathcal O}_{K}\otimes_{\mathbb Z}W(k)$-module of rank one.
%\end{remark}
%%%%%%%
%In {\bf Theorem 1.1} if $t=g$, $k={\mathbb F}_{p}$ and if $A$ admits to a lift to ${\mathbb Z}_{p}$, we will obtain a stronger result.
\begin{thm} Let ${A}$ be an abelian variety over ${\mathbb F}_{p}$ of dimension $g$ with CM by ${\mathcal O}_{K}$ and we assume that the following conditions are satisfied:
\begin{enumerate}
\item $p$ completely splits in the maximal totally real subfield $K_{0}$:
\[p={\frak P}_{1}\cdots{\frak P}_{g},\]
and that each prime ${\frak P}_{i}$ remains prime in $K$. 
\item The action of ${\mathcal O}_{K_{0}}$ is defined over ${\mathbb F}_{p}$. 
\end{enumerate}
%Let $A$ be an abelian variety over ${\mathbb F_{p}}$ of dimension $g$ with CM which lifts to ${\mathbb Z}_{p}$. Suppose that $p$ completely splits in $K_{0}$:
%\[p={\frak P}_{1}\cdots{\frak P}_{g},\]
%and that every ${\frak P}_{i}$ remains prime in $K$. 
Then $A$ is a product of supersingular elliptic curve $\{E_{i}\}_{1\leq i \leq g}$, 
\[A=E_{1}\times  \cdots \times E_{g},\]
over ${\mathbb F}_{p}$. If moreover $p$ is greater than $3$ we have $\Phi_{A,p}(t)=(t^2+p)^g$.
\end{thm}
A projective smooth curve $C$ of genus $g$ defined over a number field $F$  is called {\it a CM-curve} if the endomorphism ring of the Jacobian variety ${\rm Jac}(\overline{C})$ contains ${\mathcal O}_{K}$, where $K$ is a CM-field satisfying $[K:{\mathbb Q}]=2g$. We call a finite prime $v$ of $F$ is {\it good} if the reduction $C_{v}$ is nonsingular. A good prime $v$ is mentioned {\it ordinary} (resp. {\it supersingular}, {\it superspecial}) if so is ${\rm Jac}(C_{v})$.
{\bf Theorem 1.1}, {\bf Theorem 1.2}, {\bf Theorem 1.3} and the Hasse-Weil's formula yield the following consequence, which generalizes the Deuring's results.
\begin{thm} Let $C$ be a proper smooth curve of genus $g$ over ${\mathbb Q}$ with CM by ${\mathcal O}_{K}$ and $p$ a good prime. 
\begin{enumerate}
\item If $p$ completely splits in $K$, $C_{p}$ is ordinary.
\item Let 
\[p={\frak P}_{1}\cdots {\frak P}_{t},\]
be the prime factorization in $K_{0}$. If every ${\frak P}_{i}$ remains prime in $K$, $C_{p}$ is supersingular. If moreover $t=g$ (i.e. $p$ completely splits in $K_{0}$) it is superspecial.
\item Suppose that $t=g$ in (2) and that either the following (a) or (b) is satisfied.
\begin{enumerate}
\item The action of ${\mathcal O}_{K_{0}}$ on ${\rm Jac}(C)$ is defined over $K_{0}$.
\item The action of ${\mathcal O}_{K}$ on ${\rm Jac}(C)$ is defined over $K$.
\end{enumerate}
Then ${\rm Jac}(C_{p})$ is a product of supersingular elliptic curves over ${\mathbb F}_{p}$. If moreover $p\geq 5$, the number of ${\mathbb F}_{p^2}$-points attains the Hasse-Weil upper bound:
\[|C({\mathbb F}_{p^{2}})|=1+p^{2}+2gp.\]
\end{enumerate}
\end{thm}

The following corollaries are special cases of of this theorem. Let $C$ be the curve in {\bf Theorem 1.4}.
\begin{cor} Suppose that $K$ is a cyclotomic field ${\mathbb Q}(\zeta_{N})$ where $\zeta_{N}$ is a primitive $N$-th root of unity that satisfies $\phi(N)=2g$, where $\phi$ is the Euler function. 
\begin{enumerate}
\item If $p\equiv 1(N)$, $C_{p}$ is ordinary.
\item Suppose that there is a positive integer $h$ with $p^{h}\equiv -1(N)$.
Then $C_{p}$ is supersingular. 
\item If $p\equiv -1(N)$, ${\rm Jac}(C_{p})$ is a product of supersingular elliptic curves over ${\mathbb F}_{p}$. If moreover $p\geq 5$,
\[|C({\mathbb F}_{p^{2}})|=1+p^{2}+2gp.\]
\end{enumerate}
\end{cor}

\begin{cor} Suppose that $K={\mathbb Q}(\zeta_{M}+\zeta_{M}^{-1}, \zeta_{d})$, where $M$ is a positive integer satisfying $\phi(M)=2g$ and $d=3$ or $4$.
\begin{enumerate}
\item Assume that $p\equiv 1(d)$ and that $p\equiv \pm 1(M)$. Then $C_{p}$ is ordinary.
\item If $p\equiv -1(d)$, $C_{p}$ is supersingular. 
\item Suppose that $p\equiv -1(d)$ and that $p\equiv \pm 1(M)$. Then ${\rm Jac}(C_{p})$ is a product of supersingular elliptic curves over ${\mathbb F}_{p}$ and if moreover $p \geq 5$, 
\[|C({\mathbb F}_{p^{2}})|=1+p^{2}+2gp.\]

\end{enumerate}
\end{cor}
In the final section we will show several curves over ${\mathbb F}_{p}$ whose the number of ${\mathbb F}_{p^2}$-points attains the Hasse-Weil's upper bound.  We hope that our theorems may offer a new construction of such {\it a maximal curve}, i.e. the number of ${\mathbb F}_{p^2}$-points attains the Hasse-Weil's upper bound.\\

Let us briefly explain how the theorems will be proved. We will reduce a problem of an abelian variety to one of the $p$-divisible group. For simplicity suppose that $p$ completely splits in $K_{0}$, $p={\frak P}_{1}\cdots{\frak P}_{g}$. 
Then completion $({\mathcal O}_{K_{0}})_{{\frak P}_{i}}$ is isomorphic to ${\mathbb Z}_{p}$ and
\[{\mathcal O}_{K_{0}}\otimes_{\mathbb Z}{\mathbb Z}_{p}\simeq ({\mathcal O}_{K_{0}})_{{\frak P}_{1}}\times \cdots \times({\mathcal O}_{K_{0}})_{{\frak P}_{g}}, \quad ({\mathcal O}_{K_{0}})_{{\frak P}_{i}} \simeq {\mathbb Z}_{p}.\]
Since $\overline{A}[p^{\infty}]$ is naturally acted by ${\mathcal O}_{K_{0}}\otimes_{\mathbb Z}{\mathbb Z}_{p}$ this yields a decomposition,
\[\overline{A}[p^{\infty}]=G_{{\frak P}_{1}}\times \cdots \times G_{{\frak P}_{g}},\]
by the Diuedonn\'{e} theory and it will proved that the height of $G_{{\frak P}_{i}}$ is two for all $i$.
%Via the natural action of ${\mathcal O}_{K_{0}}\otimes_{\mathbb Z}{\mathbb Z}_{p}$ on the $p$-divisible group $\overline{A}[p^{\infty}]$ this induces a decomposition:
%\[\overline{A}[p^{\infty}]=G_{{\frak P}_{1}}\times \cdots \times G_{{\frak P}_{g}}.\]
In order to determine them we will study the action of the Frobenius and the Verschiebung on the Dieudonn\'{e} module of each factor. This will connect a type of $G_{{\frak P}_{i}}$ and a decomposing  pattern of ${\frak P}_{i}$ in $K$. In fact using the classification of $p$-divisible group we will show that $G_{{\frak P}_{i}}$ is isomorphic to $G_{1.1}$ or $G_{1.0}\times G_{0,1}$ according to whether ${\frak P}_{i}$ remains prime or completely splits in $K$, respectively ({\bf Lemma 3.1} and {\bf Lemma 3.2}). Now {\bf Theorem 1.1} and {\bf Theorem 1.2} will follow from the facts that relate the type of $\overline{A}$ and the shape of $\overline{A}[p^{\infty}]$ (\cite{Deligne}, \cite{Oort1974}, \cite{Oort1975}, \cite{Shioda}, \cite{Tate1966}, see {\bf Fact 2.2},  {\bf Fact 2.3} and {\bf Fact 2.4} below). In order to prove {\bf Theorem 1.3} it is sufficient to show that the product in {\bf Theorem 1.1} is defined over ${\mathbb F}_{p}$. We will deduce it from the corresponding decomposition of the $p$-divisible group ({\bf Proposition 4.1}) and {\bf Theorem 1.4} will be a consequence of preceding theorems. But a little care is needed to show (3) and here the assumption that $C_p$ is a reduction of a curve over ${\mathbb Q}$ is neccessary. That is we have to check that the action of ${\mathcal O}_{K_{0}}$ on ${\rm Jac}(C_{p})$ is defined over ${\mathbb F}_{p}$. Since a simple observation shows that (b) implies (a), we may assume that (a) is satisfied. Let ${\mathcal J}$ be the Neron's model of ${\rm Jac}(C)$ over ${\mathbb Z}_{p}$, which is an abelian scheme. We claim the action of ${\mathcal O}_{K_{0}}$ is defined over ${\mathbb Z}_{p}$. By the Neron's mapping property our claim is true if the action on the generic fiber ${\rm Jac}(C)\otimes_{{\mathbb Q}}{\mathbb Q}_{p}$ is defined over ${\mathbb Q}_{p}$. This will be checked by the faithful representation of ${\mathcal O}_{K_{0}}$ on the cotangent space of ${\rm Jac}(C)\otimes_{{\mathbb Q}}{\mathbb Q}_{p}$ at the origin, which is identified with $H^0(C,\Omega^1)\otimes_{\mathbb Q}{\mathbb Q}_p$. Now consider the special fiber and the rationality of the action will be obtained.\\

Let us mention precedent results that generalize Deuring's results. Let $A$ be an abelian surface with CM by the integer ring of a cyclic quartic CM-field $K$ defined over a number field and $A_v$ its reduction at a good prime $v$ over $p$. Goren\cite{Goren} has shown that if $p$ completely splits in $K_{0}$, $p={\frak P}_{1}{\frak P}_{2}$ and if each ${\frak P}_{i}$ remains prime in $K$, $\overline{A_p}$ is a product of supersingular elliptic curves. This coincides with {\bf Theorem 1.1}. But he has also proved that $\overline{A_p}$ is simple and ordinary (resp. isogeneous but not isomorphic to a product of supersingular elliptic curves) if $p$ completely splits (resp. remains prime) in $K$. Using the Kraft's diagrams, Zaytsev\cite{Zaytsev} has completely determined $p$-torsion group of an abelian variety over a finite field $k$ (${\rm Char}\,k=p$) with the dimension less than $4$. Our {\bf Theorem 1.1} and {\bf Theorem 1.2} are contained in his results if the dimension of the abelian variety is less than $4$.\\

% Also in \cite{Sugiyama} we have obtained partial results of the present paper. Namely there we have proved {\bf Theorem 1.2} and the last half of {\bf Theorem 1.1} by a different way.  Here let us make a correction of \cite{Sugiyama}. At first we have to assume that the prime $p$ should be greater than $3$. Moreover {\bf Theorem 1.3} of \cite {Sugiyama} says that an abelian variety with CM over a finite field whose $p$-divisible group is a product of ${\mathbb L}$ is isomorphic an product of ordinary elliptic curves. But this is incorrect. In fact it is known that  there is an absolutely simple abelian variety $B$ with $B[p^{\infty}]\simeq {\mathbb L}^{\dim B}$ \cite{Lenstra-Oort}. \\

\noindent{\bf Acknowledgement}. The author appreciates Professor A. Zaytsev for his interest and useful remarks. He is also grateful for the generous referee who points out several mistakes and kindly suggests many improvements. This research is partially supported by JSPS grants Kiban(C)22540068.

%We say a smooth proper curve over a number field to be with {\it CM} if so is the Jacobian. 

%\[p={\frak P}^{0}_{1}\cdots{\frak P}^{0}_{r}{\frak Q}^{0}_{1}\cdots{\frak Q}^{0}_{s},\]
%is the prime factorization in $K_{0}$ such that each ${\frak P}^{0}_{i},\,(1\leq i \leq r)$ (resp. ${\frak Q}^{0}_{j},\,(1\leq j \leq s)$) remains prime (reps. completely splits) in $K$.

\section{A review of $p$-divisible groups.}
In this section we summarize facts of $p$-divisible groups and the Dieudonn\'{e} functor which will be used later. The references are \cite{Demazure}, \cite{Grothendieck}, \cite{Manin}, \cite{M-W1971}, \cite{Oort2007} and \cite{Pink}. Throughout the section $k$ will be a field of characteristic $p$.\\

\begin{df} Let $h$ be a nonnegative integer. {\rm A $p$-divisible group $G$ of height $h(G)=h$} over $k$ is an inductive system of finite group schemes $G_{i}\rightarrow {\rm Spec}\,k,\,(i\geq 1),$ satisfying
\begin{enumerate}
\item the dimension of the coordinate ring $k[G_{i}]$ over $k$ equals to $p^{h\cdot i}$,
\item $p^{i}$ annihilates $G_{i}$,
\item there are inclusions $G_{i} \hookrightarrow G_{i+1}$ such that
\[G_{i+1}[p^{i}]=G_{i},\]
\end{enumerate}
and we denote 
\[G=\lim_{\rightarrow}G_{i}.\]
\end{df}
\begin{remark}
Let $\Gamma$ be a finite group scheme over $k$. {\rm The order} is defined to be the dimension of the coordinate ring $k[\Gamma]$ and is described by $|\Gamma|$.
% We mention the dimension of the coordinate ring of a finite group scheme $\Gamma$ over $k$ {\rm the order} and denote it by $|\Gamma|$.
% {\it the order} $|\Gamma|$ of a finite group scheme $\Gamma$ is defined to be the dimension of its coordinate ring over $k$. 
\end{remark}
In the definition,
\[X[f]={\rm Ker} [X \stackrel{f}\rightarrow X],\]
for an endomorphism $f$ of a group scheme $X$ and note that $G_{i+j}/G_{i}=G_{j}$. Let $Z$ be a connected commutative formal smooth group scheme of finite type over $k$. Then $\{Z[p^i]\}_i$ is a $p$-divisible group and the {\it dimension} is defined to be one of $Z$. It is equal to the order of the kernel of the Frobenius on $Z$(\cite{Demazure} {\bf Chapter 3}). Here are examples of $p$-divisible groups.
\begin{enumerate}
\item 
\[{\mathbb G}_{m}[p^{\infty}]:=\lim_{\rightarrow}{\mathbb G}_{m}[p^{n}],\]
whose height and dimension are $1$.

\item 
\[{\mathbb Q}_{p}/{\mathbb Z}_{p}:=\lim_{\rightarrow}{\mathbb Z}/p^{n}{\mathbb Z},\]
is the  \'{e}tale $p$-divisible group of height one:
\item Let $A$ be an abelian variety over $k$ of dimension $g$. Then
\[A[p^{\infty}]:=\lim_{\rightarrow}A[p^{n}],\]
 is a $p$-divisible group of height $2g$.
 
\end{enumerate}
Let $G=\lim_{\rightarrow}G_{i}$ be a $p$-divisible group. Then so is the collection of the Cartier dual $G^{\vee}:=\lim_{\rightarrow} G_{i}^{\vee}$ and called {\it the Serre dual}. For example ${\mathbb Q}_{p}/{\mathbb Z}_{p}$ and ${\mathbb G}_{m}[p^{\infty}]$ are Serre dual to each other. Taking the Serre dual induces an involution of the category of $p$-divisible groups. \\

From now on we assume that $k$ is perfect field of characteristic $p$.  Let $W$ be the Witt group scheme defined over $k$ of $\infty$-length. It is isomorphic to a product of infinite affine lines as a scheme 
\[\prod_{n\geq 0}{\mathbb A}^{1}={\rm Spec} k[x_{0}, x_{1}, \cdots ], \]
and we denote a point $x$ by $x=(x_{0}, x_{1}, \cdots)$. The Frobenius $F$ and the Verschiebung $V$ which are endomorphism of $W$ are defined to be
\[F((x_{0}, x_{1}, \cdots))=(x_{0}^{p}, x_{1}^{p}, \cdots), \quad V((x_{0}, x_{1}, \cdots))=(0, x_{0}, x_{1}, \cdots),\]
and since $FV=VF=p$, 
\[p((x_{0}, x_{1}, \cdots))=(0, x_{0}^{p}, x_{1}^{p}, \cdots).\]
For positive integer $n$ let $W_{n}$ denote the additive group scheme of Witt vectors of length $n$. It is isomorphic to $W/V^{n}W$ and the collection of all $\{W_{n}\}$ forms a direct system by the natural inclusions.\\

%%%%%%%%
% and we denote the kernel of $F^{m}$ on $V_{n}$ by $W_{n}^{m}$.  For an example $W_{1}$ is isomorphic to ${\mathbb A}^{1}$ and $W_{1}^{1}\simeq \alpha_{p}$, where $\alpha_{p}$ is a finite group scheme defined by $\alpha_{p}:={\rm Spec} k[x]/(x^{p})$. The collection of all $\{W^{m}_{n}\}$ forms a direct system by the following injective homomorphisms:
%\[\begin{CD}
% W^{m}_{n}@>\text{$i$}>>  W^{m+1}_{n}\\
% @V\text{$v$}VV  @V\text{$v$}VV\\
%  W^{m}_{n+1}  @>\text{$i$}>>   W^{m+1}_{n+1}.
% \end{CD}\]

 Let $W(k)$ be the ring of Witt vectors whose coefficients are in $k$. The Frobenius induces the ring homomorphism and will be denoted by $\sigma$.  Let ${\mathcal D}(k)$ be a non-commutative algebra whose coefficient ring is $W(k)$, which is generated by semi-linear operators $F$ and $V$ with the relations
\[FV=VF=p, \quad F\lambda=\lambda^{\sigma}F, \quad \lambda V=V\lambda^{\sigma},\quad \forall \lambda\in W(k).\]
Consider the torsion $W(k)$-module
\[T:=W(k)[\frac{1}{p}]/W(k).\]
Then the functor
\[N \mapsto N^{*}:={\rm Hom}_{W(k)}(N,T),\]
defines an anti-equivalence from the abelian category of finite length $W(k)$-modules to itself and 
\[N \simeq (N^{*})^{*}.\]
We define the actions of $F$ and $V$ on $N^{*}$ is defined as
\[(Fl)(n):=\sigma(l(Vn)),  \quad (Vl)(n):=\sigma^{-1}(l(Fn)), \quad l\in N^{*}, n\in N.\]
Let  $\Gamma$ be an affine unipotent group over $k$. The {\it the Dieudonn\'{e} module} of $\Gamma$ is defined to be
\[{\mathbb M}(\Gamma):=\lim_{n\to \infty}{\rm Hom}(\Gamma, W_{n}).\]
Here "${\rm Hom}$" is taken in the category of affine unipotent group schemes over $k$. It is a contra-variant functor from the category of the affine unipotent group over $k$ to that of all ${\mathcal D}(k)$-modules killed by a power of $V$ and induces an anti-equivalence between them. $\Gamma$ is algebraic (resp. finite) if and only if ${\mathbb M}(\Gamma)$ is a finitely generated ${\mathcal D}(k)$-module (resp. a $W(k)$-module of finite length). By restriction it induces an anti-equivalence between the category of finite unipotent \'{e}tale (resp. infinitesimal) groups over $k$ and that of ${\mathcal D}(k)$-modules which are $W(k)$-module of finite length, killed by a power of $V$ and on which $F$ is bijective (resp. killed by a power of $F$). Finally the Dieudonn\'{e} module of a finite infinitesimal multiplicative group $\Gamma$ is defined by
\[{\mathbb M}(\Gamma):={\mathbb M}(\Gamma^{\vee})^{*}.\]
Then the functor ${\mathbb M}$ induces an anti-equivalence between the abelian category of finite commutative group schemes over $k$ of a $p$-power order to that of left ${\mathcal D}(k)$-modules of finite length. It is known that the length of ${\mathbb M}(\Gamma)$ is equal to $\log_{p}|\Gamma|$. Here are some examples.
\begin{ex}
\begin{enumerate}
\item Let $W_{n}^{m}$ be the kernel of $F^{m}$ on $W_{n}$. Then
\[{\mathbb M}(W_{n}^{m})={\mathcal D}(k)/({\mathcal D}(k)F^{m}+{\mathcal D}(k)V^{n}).\]
\item ${\mathbb M}({\mathbb Z}/p{\mathbb Z})$ (resp. ${\mathbb M}({\mathbb G}_{m}[p])$ ) is isomorphic to $k$ with $F=1,\, V=0$ (resp. $F=0, \, V=1$).
\item Let $E$ be a supersingular elliptic curve over $k$, then
\[{\mathbb M}(E[p])={\mathcal D}(k)\otimes_{W(k)}k/{\mathcal D}(k)\otimes_{W(k)}k(F-V).\]
\end{enumerate}
\end{ex}

%%%%%%%%%
%For any finite commutative group scheme $\Gamma$ over $k$ of local type
%There is a contra-variant functor ${\mathbb M}$ called {\it the Dieudonn\'{e} functor} from an abelian category of finite group schemes over $k$ of $p$-power order to one of left ${\mathcal D}(k)$-modules of finite length. (Here {\it the order} $|\Gamma|$ of a finite group scheme $\Gamma$ is defined to be the dimension of its coordinate ring over $k$). It is known that the length of ${\mathbb M}(\Gamma)$ is equal to $\log_{p}|\Gamma|$. Due to Dieudonn\'{e}, Cartier, Barsotti and Oda, it gives an anti-equivalence of these categories (\cite{M-W1971} {\bf Theorem 4}). \\

%
%Namely for a finite group scheme $\Gamma$ of $p$-power order we have
%\[{\mathbb M}(\Gamma^{\vee})={\mathbb M}(\Gamma)^{\vee},\]
%where $\Gamma^{\vee}$ (resp. ${\mathbb M}(\Gamma)^{\vee}$) is the Cartier (resp. Pontryagin) dual. 
%
We define {\it the Dieudonn\'{e} module} of a $p$-divisible group $G=\lim_{\rightarrow}G_{i}$ as
\[ {\mathbb M}(G):=\lim_{\leftarrow}{\mathbb M}(G_{i}),\]
which is a ${\mathcal D}(k)$-module by definition. It is a free $W(k)$-module with rank $h(G)$. In this way the Dieudonn\'{e} functor ${\mathbb M}$ gives a contra-equivalence between the category of $p$-divisible groups defined over $k$ and one of ${\mathcal D}(k)$-modules that are free over $W(k)$ with finite rank. For a pair of coprime integers $(d,c)$ so that $d>0,\,c\geq 0$ the $p$-divisible group $G_{d,c}$ is defined to be 
\[G_{d,c}:={\rm Ker}[F^{c}-V^{d}: W[p^{\infty}] \to W[p^{\infty}] ],\]
where $W[p^{\infty}]:=\lim_{n\to \infty} W[p^{n}]$. Its Dieudonn\'{e} module is 
\begin{equation}{\mathbb M}(G_{d,c})={\mathcal D}(k)/{\mathcal D}(k)(F^{c}-V^{d}),\end{equation}
and the dimension and the height are $d$ and $c+d$, respectively. 
%%%%%%%%
% Moreover it is simple in the sense that any epimorphism is either an isogeny or the structure morphism $G_{d,c} \to {\rm Spec}k$. Here we call a homomorphism of $p$-divisible groups isogeny if both kernel and cockerel are finite. 
 %%%%%%%%%%%%%
 One sees that $G_{1,0}\simeq {\mathbb G}_{m}[p^{\infty}]$ and it is convenient to set $G_{0,1}:={\mathbb Q}_{p}/{\mathbb Z}_{p}$. 
% It has dimension $n$ and any epimorphism is either an isogeny or the structure morphism $G_{n,m} \to {\rm Spec}k$. Here we call homomorphism of $p$-divisible groups isogeny if both kernel and cockerel are finite. Following Manin \cite{Manin} we use the following convention. We denote ${\mathbb G}_{m}[p^{\infty}]$ (resp. ${\mathbb Q}_{p}/{\mathbb Z}_{p}$) by $G_{1,0}$ (resp. $G_{0,1}$) and the $p$-divisible group of the Dieudonn\'{e} module ${\mathcal D}(k)/{\mathcal D}(k)V^{n}$ by $G_{n,\infty}$, which is derived from the formal group of $W_{n}$. Note that its height is infinity. 
 %Temporary we assume that $k$ is an algebraic closed field.  The dimension and the height of $G_{d,c}$ is $d$ and $c+d$, respectively.  
%By definition $G_{c,d}=(G_{d,c})^{\vee}$ and $G_{d,c}$ satisfies an equation:
%\[
%G_{d,c}[F^{c+d}]=G_{d,c}[p^{d}],\quad G_{d,c}[V^{c+d}]=G_{d,c}[p^{c}].
%\]
Then $G_{c,d}=(G_{d,c})^{\vee}$ for every pair of coprime non-negative integer $(d,c)\neq(0,0)$. 
Temporally we assume that $k$ is algebraically closed. Then $G_{d,c}$ is characterized by the height $h=c+d$ and {\it the slope} $d/h=d/(c+d)$. A $p$-divisible group $G$ over $k$ is mentioned {\it simple} if any epimorphism from $G$ is either an isogeny or the structure morphism (here an isogeny is a homomorphism whose kernel and cockerel are finite). Then $G_{d,c}$ is simple and conversely any simple $p$-divisible group is isomorphic to a certain $G_{d,c}$.
%Moreover it is simple in the sense that any epimorphism whose source is $G_{d,c}$ is either an isogeny or the structure morphism $G_{d,c} \to {\rm Spec}k$. Here we call a homomorphism of $p$-divisible groups isogeny if both kernel and cockerel are finite.  Any simple $p$-divisible group $G$ is isomorphic to a certain $G_{d,c}$. 
Note that a simple $p$-divisible group $G$ is isomorphic to $G_{d,c}$ if and only if there is a pair of non-negative integer $(m,n)\neq (0,0)$ so that
\begin{equation}
G[F^{m}]=G[V^{n}], \quad \frac{n}{m+n}=\frac{d}{c+d}.
\end{equation}
In fact (2) shows that $G\simeq G_{d,c}$ if and only if
 \[{\mathbb M}(G)/V^{n}{\mathbb M}(G)={\mathbb M}(G)/F^{m}{\mathbb M}(G),\]
 for a pair of non-negative integer $(m,n)\neq (0,0)$  satisfying $d/(c+d)=n/(m+n)$ and this is equivalent to (3). Here are $p$-divisible groups whose height is less than $3$ (\cite{Demazure}, p.93):
\begin{enumerate}
\item $h(G)=0$ iff $G=0$. 
\item If $h(G)=1$, then $G=G_{1.0}$, or $G_{0.1}$.
\item If $h(G)=2$, $G$ is the one of the followings,
\[G_{1.0}^{2}, \quad G_{1.0}\times G_{0.1}, \quad G_{0.1}^{2}, \quad G_{1.1}.\]
\end{enumerate}
Set ${\mathbb L}= G_{1.0}\times G_{0.1}$. Then $G_{1.1}$ (resp. ${\mathbb L}$) is isomorphic to the $p$-divisible group of a supersingular (resp. an ordinary) elliptic curve. Let $X$ and $Y$ be $p$-divisible groups. If there is an isogeny between them we say that  they are {\it isogeneous} and describe as $X\sim Y$. This notion defines an equivalence relation on the set of $p$-divisible groups.
\begin{fact}( \cite{Demazure}, p.85 or \cite{Manin}, p.35) Let $G$ be a $p$-divisible group over an algebraically closed field $k$. Then there is an isogeny:
\[G \sim G_{1,0}^{f}\times G_{0,1}^{f^{\prime}}\times \prod_{i}G_{d_i,c_i},\]
where $\{d_i,c_i \}_{i}$ are pairs of positive coprime integers.
\end{fact}

Let $A$ be an abelian variety over $k$ of dimension $g$ and set $\overline{A}:=A\otimes_{k}\overline{k}$. We define {\it $p$-rank} $f(A)=f$ as an integer such that $A[p](\overline{k})\simeq ({\mathbb Z}/p{\mathbb Z})^{f}$, which equals to $\dim_{{\mathbb F}_{p}}{\rm Hom}({\mathbb G}_{m}[p], \overline{A}[p])$. Let $\alpha_{p}$ be a finite group scheme defined by $\alpha_{p}:={\rm Spec}\,\overline{k}[x]/(x^{p})$. The {\it $a$-number} is defined to be $a(A):=\dim_{\overline{k}}{\rm Hom}(\alpha_{p}, \overline{A}[p])$.
%Then we set $a(A)=\dim_{\overline{k}}{\rm Hom}(\alpha_{p}, \overline{A}[p])$ and call it the {\it $a$-number}. 
It is known that $0\leq f(A) \leq g$ and that $0\leq a(A)+f(A) \leq g$. We say $A$ {\it ordinary} if $f(A)=g$. 
\begin{fact}(\cite{Deligne}) Let $A$ be an abelian variety over $k$ of dimension $g$. Then the following are equivalent:
\begin{enumerate}
\item $A$ is ordinary,
\item $\overline{A}[p^{\infty}]\sim {\mathbb L}^{g}$.
\end{enumerate} 
\end{fact}
This is well-known if $g=1$. In fact let $E$ be an elliptic curve defined over an algebraic closed field of characteristic $p$. Since there is an exact sequence
\begin{equation}0\to \alpha_{p} \to G_{1.1}[p] \to \alpha_{p} \to 0,\end{equation}
the previous classification of $p$-divisible group of height $2$ shows that $a(E)=1$ if and only $E$ is supersingular. This observation is generalized to a higher dimensional abelian varieties. The following fact is due to Deligne, Oort, Shioda and Tate.

\begin{fact}(\cite{Oort1974}, \cite{Shioda}, \cite{Tate1966}) Let $A$ an abelian variety over an algebraic closed field $k$ with dimension $g$. If
\[A[p^{\infty}]\sim G_{1.1}^{g},\]
$A$ is supersingular. 
\end{fact}

\begin{fact}(\cite{Oort1975} {\bf Theorem 2}) Let $A$ an abelian variety over an algebraic closed field $k$ of with dimension $g$. If
\[a(A)=g,\]
$A$ is superspecial. 
\end{fact}

Let $A$ be an abelian variety over $k$ of dimension $g$. We denote the Dieudonn\'{e} module of $A[p^{\infty}]$ by $T_{p}(A)$, which is a free $W(k)$-module of rank $2g$.

%For an abelian variety $A$ over $k$ of dimension $g$, following Milne and Waterhaus\cite{M-W1971}, we denote ${\mathbb M}(A[p^{\infty}])$ by $T_{p}(A)$. 
\begin{fact} (\cite{M-W1971} {\bf Theorem 6}) Let $A$ and $B$ are abelian varieties over a finite field $k$. Then 
 \[{\rm Hom}_{k}(A,B)\otimes_{\mathbb Z} {\mathbb Z}_{p} \simeq {\rm Hom}_{{\mathcal D}(k)}(T_{p}(B), T_{p}(A)).\]
 \end{fact}

% For an abelian variety $A$ over $k$ of dimension $g$ we denote ${\mathbb M}(A[p^{\infty}])$ by $T_{p}(A)$. Since the height of $A[p^{\infty}]$ is $2g$ $T_{p}(A)$ is a free $W(k)$-module of rank $2g$.
 
% \begin{fact}(\cite{Grothendieck}  IV \$7 and V \$5) Let $k$ be a perfect field of characteristic $p$  and $A$ the special fiber of an abelian scheme ${\mathcal A}$ over $W(k)$. Then $T_{p}(A)$ is isomorphic to the first de Rham cohomology group $H^{1}_{{\rm DR}}({\mathcal A}/W(k))$ as ${\mathcal D}(k)$-modules.
% \end{fact}
%%%%%%%%%%%%
\section{The $p$-divisible group of an abelian variety with CM}
Let $k$ be a finite field of characteristic $p$ with $[k:{\mathbb F}_{p}]=r$ and $A$ an abelian variety over $k$ of dimension $g$ endowed with CM by ${\mathcal O}_{K}$. We fix an imbedding of the CM-field $K$ into ${\mathbb C}$ and denote the restriction of the complex conjugation to $K$ by "$\prime$". We assume that ${\rm End}_{k}(A)={\rm End}_{\overline{k}}(\overline{A})$ and denote it by $R$. We first note that $F^{r}$ and $V^{r}$ are contained in ${\mathcal O}_{K}$. In fact here is a proof after \cite{Zaytsev} {\bf Lemma 3.2}. Set $\pi_{F}:=F^{r}$ and $\pi_{V}:=V^{r}$. The assumption is equivalent to that $\pi_{F}$ and $\pi_{V}$ are contained in the center $C$ of $R\otimes_{\mathbb Z}{\mathbb Q}$. Since the commutant of $K$ in $R\otimes_{\mathbb Z}{\mathbb Q}$ is itself (\cite{Serre-Tate} \$4 {\bf Corollary 1}),  $C$ is contained in $K$ and $\pi_{F}, \,\pi_{V} \in {\mathcal O}_{K}=K\cap R$. It is known that $\pi_{F}$ is a $q$-Weil number ($q=p^{r}$) and $\pi_{F}\cdot \pi_{F}^{\prime}=p^{r}$ (\cite{Oort2007} {\bf Theorem 3.2} and {\bf Proposition 2.2}). On the other hand since $FV=VF=p$, $\pi_{F}\cdot \pi_{V}=F^{r}V^{r}=p^{r} $. Therefore $(\pi_{F}^{\prime}-\pi_{V})\pi_{F}=0$ and because $\pi_{F}$ is surjective we conclude
\[\pi_{F}^{\prime}=\pi_{V}.\]
For simplicity we will denote $\pi_{F}$ and $\pi_{V}$ by $\pi$ and $\pi^{\prime}$, respectively. 
%Since the height of $A[p^{\infty}]$ is $2g$, $T_p(A)$ is a free $W(k)$-module of rank $2g$. Moreover by {\bf Fact 2.5} it has a faithful action of ${\mathcal O}_K\otimes_{\mathbb Z}W(k)$.
%%%%%%%%%%%%%%%%%%%%
%\begin{lm}$T_{p}(A)$ is a free $W(k)$-module with rank $2g=[K:{\mathbb Q}]$ and has a faithful action of ${\mathcal O}_K\otimes_{\mathbb Z}{\mathbb Z}_p$.
%\end{lm}
%\noindent {\bf Proof.}
%Let ${\mathcal D}_{r}(k)$ be a subalgebra of ${\mathcal D}(k)$ which is generated by $\pi$ and $\pi^{\prime}$ over $W(k)$. Notice that $\sigma^{r}$ is the identity on $W(k)$ and ${\rm End}_{{\mathcal D}_{r}(k)}T_{p}(A)$ contains $W(k)$-linear extension of ${\rm End}_{{\mathcal D}(k)}T_{p}(A)$. Therefore the following monomorphisms are obtained:
%\[({\rm End}_{{\mathcal D}(k)}T_{p}(A))\otimes_{{\mathbb Z}_{p}}W(k) \hookrightarrow {\rm End}_{{\mathcal D}_{r}(k)}T_{p}(A)\hookrightarrow {\rm End}_{W(k)}T_{p}(A).\]
%Since ${\mathcal O}_{K}$ is contained in ${\rm End}_{k}A$, {\bf Fact 2.5} yields an injective homomorphsim
%\[ {\mathcal O}_{K}\otimes_{\mathbb Z}W(k)\hookrightarrow {\rm End}_{W(k)}T_{p}(A).\]
%This shows that $T_{p}(A)$ is a free ${\mathcal O}_{K}\otimes_{\mathbb Z}W(k)$-module whose rank is one because
%\[{\rm rank}_{W(k)}T_{p}(A)={\rm rank}_{W(k)}{\mathcal O}_{K}\otimes_{\mathbb Z}W(k)=2g.\]
%\begin{flushright}
%$\Box$
%\end{flushright}
%Throughout the following we will identify $T_{p}(A)$ with ${\mathcal O}_{K}\otimes_{\mathbb Z}W(k)$ as modules by {\bf Lemma 3.1}.  
As we have seen in (3) a simple $p$-divisible group is characterized by its slope. Since $\pi=F^{r}$ and $\pi^{\prime}=V^{r}$ are contained in ${\mathcal O}_{K}$,  we will obtain an information of the Dieudonn\'{e} module of a simple component of $A[p^{\infty}]$ from the behavior of $\{\pi,\,\pi^{\prime}\}$ in ${\mathcal O}_{K}\otimes_{\mathbb Z}{\mathbb Z}_p$ with help of {\bf Fact 2.5}. This is our strategy.\\

%%%%%%%%%%%%%%%%%%%%
%We assume that the pair $(A,{\mathcal O}_{K})$ has a lift $({\mathcal A}, {\mathcal O}_{K})$ to $W(k)$ and that the action is defined over $W(k)$. Since the action of ${\mathcal O}_{K}$ on ${\mathcal A}$ is faithful $H^{1}_{{\rm DR}}({\mathcal A}/W(k))$ is a free ${\mathcal O}_{K}\otimes_{\mathbb Z}W(k)$-module.  The facts that ${\rm rank}_{W(k)}H^{1}_{{\rm DR}}({\mathcal A}/W(k))=2g$ and {\bf Fact 2.5} imply $T_{p}(A)$ is a free ${\mathcal O}_{K}\otimes_{\mathbb Z}W(k)$-module of rank one.
%\begin{remark} It seems that $T_{p}(A)$ may be a free ${\mathcal O}_{K}\otimes_{\mathbb Z}W(k)$-module of rank one without the assumption of a lift.
%\end{remark}
%Note that
%$T_{p}(A)$ is a free ${\mathcal O}_{K}\otimes_{\mathbb Z}W(k)$ of rank one. In fact ${\mathcal O}_{K}\otimes_{\mathbb Z}{\mathbb Z}_{p}$ is contained in ${\rm End}_{k}(A)\otimes_{\mathbb Z} {\mathbb Z}_{p}$ and {\bf Fact 2.5} implies that $T_{p}(A)$ is a torsion free ${\mathcal O}_{K}\otimes_{\mathbb Z}{\mathbb Z}_{p}$-module. Since the height of $A[p^{\infty}]$ is $2g$ $T_{p}(A)$ is a free $W(k)$-dodule of rank $2g$, which is equal to ${\rm rank}_{{\mathbb Z}_{p}}{\mathcal O}_{K}\otimes_{\mathbb Z}{\mathbb Z}_{p}$ have the same rank $2g$ over ${\mathbb Z}_{p}$, $T_{p}(\overline{A})$ is a free ${\mathcal O}_{K}\otimes_{\mathbb Z}{\mathbb Z}_{p}$ of rank one. 
%
%  Then $H^{1}_{DR}({\mathcal A}/W(\overline{\mathbb F}_{p}))$ is isomorphic to ${\mathcal O}_{K}\otimes_{\mathbb Z}W(\overline{\mathbb F}_{p})$. 
%  
 We assume that $p$ is unramified in $K$ and let 
\[p={\frak P}_{1}\cdots {\frak P}_{t},\]
be the prime factorization in $K_{0}$ so that ${\frak P}_{i}$ remains prime (resp. completely splits) in $K$ for $1\leq i \leq s$ (resp. $s+1\leq i \leq t$). Thus 
\begin{equation}p=\tilde{\frak P}_{1}\cdots \tilde{\frak P}_{s}\tilde{{\frak P}}_{s+1}\tilde{{\frak P}}^{\prime}_{s+1}\cdots \tilde{{\frak P}}_{t}\tilde{{\frak P}}^{\prime}_{t},\end{equation}
in $K$. We set ${\mathcal P}_{inert}:=\{{\frak P}_{1},\cdots, {\frak P}_{s}\}$ and ${\mathcal P}_{inert}^{f=1}:=\{{\frak P} \in {\mathcal P}_{inert}\,|\, f({\frak P}/p)=1\}$, where $f({\frak P}/p)$ is the inertia degree. Similarly ${\mathcal P}_{split}:=\{{\frak P}_{s+1},\cdots, {\frak P}_{t}\}$ and ${\mathcal P}_{split}^{f=1}:=\{{\frak P} \in {\mathcal P}_{split}\,|\, f({\frak P}/p)=1\}$. By (5),
\[{\mathcal O}_{K}\otimes_{\mathbb Z}{\mathbb Z}_{p} \simeq W({\mathbb F}_{\tilde{\frak P}_{1}})\times \cdots \times W({\mathbb F}_{\tilde{\frak P}_{s}})\times \{W({\mathbb F}_{\tilde{{\frak P}}_{s+1}})\times W({\mathbb F}_{\tilde{{\frak P}}^{\prime}_{s+1}})\} \times \cdots \times \{W({\mathbb F}_{\tilde{{\frak P}}_{t}})\times W({\mathbb F}_{\tilde{{\frak P}}^{\prime}_{t}})\},\]
where ${\mathbb F}_{\tilde{\frak P}}$ is the residue field. Using this we define $e_{i}$ to be an idempotent of ${\mathcal O}_{K}\otimes_{\mathbb Z}{\mathbb Z}_{p}$ which corresponds $(0,\cdots,0,1,0,\cdots,0)$ in RHS, where $"1"$ is placed at the $i$-th from the left. Remember that since the height of $A[p^{\infty}]$ is $2g$, $T_p(A)$ is a free $W(k)$-module of rank $2g$. Moreover by {\bf Fact 2.5} it has a faithful action of ${\mathcal O}_K\otimes_{\mathbb Z}W(k)$. The following lemma is a consequence of these facts.
\begin{lm} $e_{i}T_p(A)$ is a ${\mathcal D}(k)$-module that is free over $W(k)$ with rank $f(\tilde{{\frak P}}_i/p)$. 
\end{lm}
{\bf Proof.} Note that, by {\bf Fact 2.5}, $e_{i}T_p(A)$ is a ${\mathcal D}(k)$-module which is free over $W(k)$  and we only have to identify the rank. By the decomposition we find that $e_i({\mathcal O}_{K}\otimes_{\mathbb Z}W(k))\simeq W({\mathbb F}_{\tilde{\frak P}_{i}})\otimes_{{\mathbb Z}_p}W(k)$, which is a product of complete discrete valuation rings
\[e_i({\mathcal O}_{K}\otimes_{\mathbb Z}W(k))\simeq R_i^{(1)}\times \cdots 
\times R_i^{(\nu(i))},\]
so that every component is free and has a finite rank over $W(k)$. Let $\epsilon_j$ be an element of  $e_i({\mathcal O}_{K}\otimes_{\mathbb Z}W(k))$ that corresponds to $(0,\cdots,0,1,0,\cdots,0)$ in RHS as before and set
\[M_i^{(j)}=\epsilon_j(e_{i}T_p(A)).\]
It is a non-zero free $R_i^{(j)}$-module because the action of ${\mathcal O}_{K}\otimes_{\mathbb Z}W(k)$ on $T_p(A)$ is faithful. Let $\mu_i(j)$ be its rank and 
\[T_p(A)\simeq \oplus_i \oplus _j (R_i^{(j)})^{\oplus \mu_i(j)}, \quad \mu_i(j)\geq 1.\]
Since ${\mathcal O}_{K}\otimes_{\mathbb Z}W(k)\simeq \prod_i \prod_j R_i^{(j)}$ and since ${\rm rank}_{W(k)}T_p(A)={\rm rank}_{W(k)}{\mathcal O}_{K}\otimes_{\mathbb Z}W(k)=2g$ we find that $\mu_i(j)=1$ for all $i$ and $j$. Therefore
\[{\rm rank}_{W(k)}e_{i}T_p(A)=\sum_{j}{\rm rank}_{W(k)}R_i^{(j)}={\rm rank}_{W(k)}W({\mathbb F}_{\tilde{\frak P}_{i}})\otimes_{{\mathbb Z}_p}W(k)=f(\tilde{{\frak P}}_i/p).\]
\begin{flushright}
$\Box$
\end{flushright} 
%{\bf Proof.} We only have to identify the rank since the other statements are clear from the construction and by {\bf Fact 2.5}. By the decomposition we find $e_i({\mathcal O}_{K}\otimes_{\mathbb Z}W(k))\simeq W({\mathbb F}_{\tilde{\frak P}_{i}})\otimes_{{\mathbb Z}_p}W(k)$, which is a product of complete discrete valuation rings
%\[e_i({\mathcal O}_{K}\otimes_{\mathbb Z}W(k))\simeq R_i^{(1)}\times \cdots 
%\times R_i^{(\nu(i))},\]
%so that every component is free and finite over $W(k)$. Let $\epsilon_j$ be an element of  $e_i({\mathcal O}_{K}\otimes_{\mathbb Z}W(k))$ that corresponds to the projector $(0,\cdots,0,1,0,\cdots,0)$ to the $j$-th factor in the right hand side and set
%\[M_i^{(j)}=\epsilon_j(e_{i}T_p(A)).\]
%It is a non-zero free $R_i^{(j)}$-module because the action of ${\mathcal O}_{K}\otimes_{\mathbb Z}W(k)$ on $T_p(A)$ is faithful. Let $\mu_i(j)$ be its rank and 
%\[T_p(A)\simeq \oplus_i \oplus _j (R_i^{(j)})^{\oplus \mu_i(j)}, \quad \mu_i(j)\geq 1.\]
%Since ${\mathcal O}_{K}\otimes_{\mathbb Z}W(k)=\prod_i \prod_j R_i^{(j)}$ and since ${\rm rank}_{W(k)}T_p(A)={\rm rank}_{W(k)}{\mathcal O}_{K}\otimes_{\mathbb Z}W(k)$ we find that $\mu_i(j)=1$ for all $i$ and $j$. Therefore
%\[{\rm rank}_{W(k)}e_{i}T_p(A)=\sum_{j}{\rm rank}_{W(k)}R_i^{(j)}={\rm rank}_{W(k)}W({\mathbb F}_{\tilde{\frak P}_{i}})\otimes_{{\mathbb Z}_p}W(k)=f(\tilde{{\frak P}}_i/p).\]
%\begin{flushright}
%$\Box$
%\end{flushright} 
Since the Dieudonn\'{e} functor ${\mathbb M}$ gives the anti-equivalence, there is a $p$-divisible subgroup $G_i$ of $A[p^{\infty}]$ such that
\[{\mathbb M}(G_{i})=e_{i}T_{p}({A}).\]
%%%%%%
Set $\Gamma_{j}:=G_{s+(2j-1)}\times G_{s+2j}$ ($1\leq j \leq t-s$) and 
\begin{equation}A[p^{\infty}]=G_{1}\times \cdots \times G_{s}\times \Gamma_{1}\times \cdots \times \Gamma_{t-s}.\end{equation}
%and 
%\[{\mathbb M}(G_{i})=e_{i}T_{p}({A}),\]
%which is a free $e_{i}({\mathcal O}_{K}\otimes_{\mathbb Z}W(k))$ of rank one.
%which is isomorphic to $W(\overline{\mathbb F}_{p})^{2f({\frak P}_{i}/p)}$ (resp. $W(\overline{\mathbb F}_{p})^{f({\frak P}_{i}/p)}$) if $1\leq i \leq s$   (resp. $s+1\leq i \leq s+2(t-s)$). 
For a prime factor $\tilde{\frak P}$ of $p$ in $K$ we denote the corresponding factor of $A[p^{\infty}]$ by $G_{\tilde{\frak P}}$. For a prime ${\frak P}$ of $K_{0}$ dividing $p$ we define $p$-divisible subgroup $G_{\frak P}$ of ${A}[p^{\infty}]$ as follows. If ${\frak P}$ is contained in ${\mathcal P}_{inert}$ define $G_{\frak P}:=G_{\tilde{\frak P}}$ where $\tilde{\frak P}$ is the unique prime of $K$ over ${\frak P}$. 
%Let $e_{\tilde{\frak P}}$ be the corresponding idempotent. Since 
%\[e_{\tilde{\frak P}}({\mathcal O}_{K}\otimes_{\mathbb Z}W(k))\simeq W({\mathbb F}_{\tilde{\frak P}})\otimes_{{\mathbb Z}_{p}}W(k),\]
%and 
Since ${\mathbb F}_{\tilde{\frak P}}\simeq {\mathbb F}_{p^{2f(\frak{P}/p)}}$  we see that by {\bf Lemma 3.1}
\begin{equation}h(G_{\frak P})={\rm rank}_{W(k)}{\mathbb M}(G_{\tilde{\frak P}})=2f({\frak P}/p).\end{equation}
On the other hand if it splits: ${\frak P}=\tilde{\frak P}\times \tilde{\frak P}^{\prime}$,  we define $G_{\frak P}:=G_{\tilde{\frak P}}\times G_{\tilde{\frak P}^{\prime}}$.  Since ${\mathbb F}_{\tilde{\frak P}}\simeq {\mathbb F}_{\tilde{\frak P}^{\prime}} \simeq {\mathbb F}_{p^{f(\frak{P}/p)}}$, a similar observation shows 
\begin{equation}h(G_{\tilde{\frak P}})=h(G_{\tilde{\frak P}^{\prime}})=f({\frak P}/p), \quad h(G_{\frak P})=h(G_{\tilde{\frak P}})+h(G_{\tilde{\frak P}^{\prime}})=2f({\frak P}/p).\end{equation}
%%%%%%%%%%%%%
%Since $h(G_{\frak P})={\rm rank}_{W(\overline{\mathbb F}_{p})}{\mathbb M}(G_{\frak P})$, $h(G_{\tilde{\frak P}})={\rm rank}_{W(\overline{\mathbb F}_{p})}{\mathbb M}(G_{\tilde{\frak P}})$
%and $h(G_{\tilde{\frak P}^{\prime}})={\rm rank}_{W(\overline{\mathbb F}_{p})}{\mathbb M}(G_{\tilde{\frak P}^{\prime}})$ we see that 
%\begin{equation}h(G_{\frak P})=2f({\frak P}/p), \quad h(G_{\tilde{\frak P}})=h(G_{\tilde{\frak P}^{\prime}})=f({\frak P}/p).\end{equation}
%

\begin{lm} For ${\frak P}\in {\mathcal P}_{split}^{f=1}$, $\overline{G}_{\frak P}$ is isomorphic to ${\mathbb L}$.
\end{lm}
\noindent {\bf Proof.} By (8), $h(G_{\tilde{\frak P}})=h(G_{\tilde{\frak P}^{\prime}})=1$ and the classification of $p$-divisible groups shows that $\overline{G}_{\tilde{\frak P}}$ or $\overline{G}_{\tilde{\frak P}^{\prime}}$ is one of $\{G_{1.0}, G_{0.1}\}$. Remember that $G_{1.0}$ (resp. $G_{0.1}$) is characterized by the fact $V$ (resp. $F$) is an isomorphism on ${\mathbb M}(G_{1.0})$ (resp. ${\mathbb M}(G_{0.1})$).
%%%%%%%%%%%%%%%
%Set $r:=[k:{\mathbb F}_{p}]$.  Since the action of ${\mathcal O}_{K}$ is defined over $k$, $\pi:=F^{r}$ commutes with the action. Since ${\mathcal O}_{K}$ is a maximal commutative subalgebra of ${\rm End}_{k}(A)$ $\pi$ is contained in it. 
Let
\[(F^{r})=(\pi)=\tilde{\frak P}^{a}(\tilde{\frak P}^{\prime})^{a^{\prime}}\delta, \quad (\delta,\tilde{\frak P})=(\delta,\tilde{\frak P}^{\prime})=1,\]
be the factorization. Then 
 %%%%%%%%%%%%%%%
%Take a positive integer $r$ so that $\pi=F^{r}$ is contained in ${\mathcal O}_{K}$ and let
%\[(\pi)=\tilde{\frak P}^{a}(\tilde{\frak P}^{\prime})^{a^{\prime}}\delta, \quad (\delta,\tilde{\frak P})=(\delta,\tilde{\frak P}^{\prime})=1,\]
%be the factorization. Then 
%%%%%%%%%%%%%%%%
\[(V^{r})=(\pi^{\prime})=\tilde{\frak P}^{a^{\prime}}(\tilde{\frak P}^{\prime})^{a}\delta^{\prime},\]
and 
\[(p^{r})=(F^{r}V^{r})=(\pi \pi^{\prime})=(\tilde{\frak P}\tilde{\frak P}^{\prime})^{a+a^{\prime}}\delta\delta^{\prime}.\] 
By (5), $r=a+a^{\prime}$. Suppose that $\overline{G}_{\tilde{\frak P}}=G_{1.0}$. This implies that $V$ is an isomorphism on ${\mathbb M}(G_{\tilde{\frak P}})$ and so is $\pi^{\prime}$ (note that the slope of a $p$-divisible group is invariant under a base change). Since the action of  ${\mathcal O}_{K}$ on ${\mathbb M}(G_{\tilde{\frak P}})$ (resp. ${\mathbb M}(G_{\tilde{\frak P}^{\prime}})$) factors through an imbedding
 ${\mathcal O}_{K} \hookrightarrow {\mathcal O}_{K,\tilde{\frak P}}\simeq {\mathbb Z}_{p}$ (resp. ${\mathcal O}_{K} \hookrightarrow {\mathcal O}_{K,\tilde{{\frak P}}^{\prime}}\simeq {\mathbb Z}_{p}$), $\pi^{\prime}$ should be a unit in ${\mathcal O}_{K,\tilde{\frak P}}$ and $a^{\prime}=0$. Therefore $(V^{r})=(\tilde{\frak P}^{\prime})^{r}\delta^{\prime}$ and $(F^{r})=\tilde{\frak P}^{r}\delta$, which implies that $F^{r}$ is an isomorphism on ${\mathbb M}(G_{\tilde{\frak P}^{\prime}})$. Hence $\overline{G}_{\tilde{\frak P}^{\prime}}\simeq G_{0.1}$ and $\overline{G}_{\frak P}\simeq G_{1.0} \times G_{0.1}={\mathbb L}$. In the case of  $\overline{G}_{\tilde{\frak P}}=G_{0.1}$ the proof is similar.
\begin{flushright}
$\Box$
\end{flushright}

\begin{lm} For ${\frak P}\in {\mathcal P}_{inert}$, $\overline{G}_{\frak P}$ is isogeneous to
$G_{1.1}^{f({\frak P}/p)}$.
\end{lm}
\noindent {\bf Proof.} By {\bf Fact 2.1} $\overline{G}_{\frak P}$ is isogeneous to $\prod G_{d.c}$.
Using (3) we will show all simple factors are isomorphic to $G_{1.1}$. Let 
\[(F^{r})=(\pi)={\frak P}^{a}\delta, \quad ({\frak P}, \delta)=1,\]
be the factorization. Since ${\frak P}^{\prime}={\frak P}$, 
\[(V^{r})=(\pi^{\prime})={\frak P}^{a}\delta^{\prime}, \quad ({\frak P}, \delta^{\prime})=1,\]
and 
\begin{equation}\overline{G}_{\frak P}[F^{r}]=\overline{G}_{\frak P}[V^{r}].\end{equation}
Let $G$ be a simple factor of $\overline{G}_{\frak P}$. Then (9) shows
\[G[F^{r}]=G[V^{r}],\]
and $G=G_{1.1}$ by (3). Thus $\overline{G}_{\frak P}$ is isogeneous to $G_{1,1}^{f}$ and 
\[h(G_{\frak P})=2f,\]
which implies the claim by (7).
%%%%%%%%%%%
%\[(p^{r})=(F^{r}V^{r})={\frak P}^{2a}\delta\delta^{\prime}.\]
%By (6) $r=2a$ and $(F^{2a})={\frak P}^{a}\delta$, which yields
%\begin{equation}G_{\frak P}[F^{2a}]=G_{\frak P}[p^{a}].\end{equation}
%
\begin{flushright}
$\Box$
\end{flushright}
Set ${\mathcal P}={\mathcal P}_{inert}\sqcup{\mathcal P}_{split}=\{{\frak P}_{1},\cdots, {\frak P}_{t}\}$. 
\begin{prop}
\begin{enumerate}
\item Suppose ${\mathcal P}={\mathcal P}_{inert}$. Then $A$ is supersingular.
\item The $a$-number of $A$ is greater than or equal to $|{\mathcal P}_{inert}^{f=1}|$. 
\item Suppose that
\[{\mathcal P}={\mathcal P}_{inert}^{f=1}\cup {\mathcal P}_{split}^{f=1}.\]
Then
\[a(A)=|{\mathcal P}_{inert}^{f=1}|,\quad f(A)=|{\mathcal P}_{split}^{f=1}|.\]
In particular if ${\mathcal P}={\mathcal P}_{inert}^{f=1}$ (resp.  ${\mathcal P}={\mathcal P}_{split}^{f=1}$) $A$ is superspecial (resp. ordinary).
\end{enumerate}
\end{prop}
\noindent {\bf Proof.} (1) is an immediate consequence of {\bf Lemma 3.3} and {\bf Fact 2.3}. If necessary arranging the indices, ${\mathcal P}_{inert}^{f=1}=\{{\frak P}_{1},\cdots, {\frak P}_{r}\}$ ($r \leq s$) and we consider a subgroup $G_{{\frak P}_{1}}\times \cdots \times G_{{\frak P}_{r}}$ of $A[p^{\infty}]$. By (7) and  {\bf Lemma 3.3} $\overline{G}_{{\frak P}_{i}}$ is isomorphic to $G_{1.1}$ for $1\leq \forall i \leq r$. Therefore $\overline{A}[p]$ contains $G_{1.1}[p]^{s}$ and by (4),
\[a(A)=\dim_{\overline{\mathbb F}_{p}}{\rm Hom}(\alpha_{p}, \overline{A}[p]) \geq \dim_{\overline{\mathbb F}_{p}}{\rm Hom}(\alpha_{p}, G_{1.1}[p]^{r})=r.\]
Finally we prove (3). Set $s=|{\mathcal P}_{inert}^{f=1}|$ and $g-s=|{\mathcal P}_{split}^{f=1}|$. {\bf Lemma 3.2} and {\bf Lemma 3.3} show that the product (6) becomes
\[\overline{A}[p^{\infty}]\simeq G_{1.1}^{s}\times {\mathbb L}^{g-s},\]
which implies $a(A)=|{\mathcal P}_{inert}^{f=1}|$ and $f(A)=|{\mathcal P}_{split}^{f=1}|$. The last statement follows from {\bf Fact 2.2} and {\bf Fact 2.4}, respectively.

\begin{flushright}
$\Box$
\end{flushright}
Now {\bf Theorem 1.1} and {\bf Theorem 1.2} are direct consequences of {\bf Proposition 3.1}.
\begin{cor} Suppose that $K$ is a Galois extension of ${\mathbb Q}$.
\begin{enumerate}
\item If ${\mathcal P}_{inert}$ is not empty, $A$ is supersingular.
\item Suppose that $p$ completely splits in $K_{0}$. If it completely splits even in $K$, $A$ is ordinary and otherwise $A$ is superspecial.
\end{enumerate}
\end{cor} 
\section {Rationality}
Let ${A}$ be an abelian variety over ${\mathbb F}_{p}$ which satisfies the assumption of {\bf Theorem 1.3}. 
%of dimension $g$ with CM by ${\mathcal O}_{K}$. We assume that the following conditions are satisfied:
%\begin{enumerate}
%\item $p$ completely splits in the maximal totally real subfield $K_{0}$:
%\[p={\frak P}_{1}\cdots{\frak P}_{g},\]
%and that every prime ${\frak P}_{i}$ remains prime in $K$. 
%\item The action of ${\mathcal O}_{K_{0}}$ is defined over ${\mathbb F}_{p}$. 
%\end{enumerate}
The completion $({\mathcal O}_{K_0})_{\frak{P}_{i}}$of ${\mathcal O}_{K_0}$ at ${\frak P}_{i}$ is isomorphic to ${\mathbb Z}_{p}$ and
\begin{equation}{\mathcal O}_{K_0}\otimes_{\mathbb Z}{\mathbb Z}_{p}\simeq {\mathbb Z}_{p}\times \cdots \times {\mathbb Z}_{p}, \quad \alpha=(\alpha_{1},\cdots,\alpha_{g}).
\end{equation}
Let $e^{0}_{i}$ be the idempotent in
${\mathcal O}_{K_0}\otimes_{\mathbb Z}{\mathbb Z}_{p}$ which corresponds
$(0,\cdots,0,1,0,\cdots, 0)$ in RHS of (10) (where $"1"$ is placed at the $i$-th from the left, as before) and $G^{0}_{i}$ a $p$-divisible subgroup of $A[p^{\infty}]$ such that ${\mathbb M}(G^{0}_{i})=e^{0}_{i}T_p(A)$.
Then
\begin{equation}A[p^{\infty}]=G^{0}_{1}\times \cdots \times G^{0}_{g},\end{equation}
and since the action of ${\mathcal O}_{K_0}\otimes_{\mathbb Z}{\mathbb Z}_{p}$ on $A[p^{\infty}]$ defined over ${\mathbb F}_{p}$ so is the product. The following lemma is clear from {\bf Lemma 3.3} (see also the proof of {\bf Proposition 3.1}).
\begin{lm}
$\overline{G}^{0}_{i}$ is isomorphic to $G_{1,1}$ for all $i$.
\end{lm}
%\noindent {\bf Proof.} We have
%\[\overline{G}^{0}_{i}=\overline{e_{i}^{0}A[p^{\infty}]}=e^{0}_{i}(\overline{A}[p^{\infty}]),\]
%and this is isomorphic to $G_{1,1}$ as we have shown in {\bf Proposition 3.1}.
%Since
%\[H^{1}_{DR}({\mathcal A}/{\mathbb Z}_{p}) \simeq {\mathcal O}_{K}\otimes_{{\mathbb Z}}{\mathbb Z}_{p} \simeq {\mathcal O}_{K}\otimes_{{\mathcal O}_{K_{0}}}({\mathcal O}_{K_{0}}\otimes_{{\mathbb Z}}{\mathbb Z}_{p})\simeq W({\mathbb F}_{p^2})^{g},\]
%and since ${\mathbb M}(A[p^{\infty}])=H^{1}_{DR}({\mathcal A}/{\mathbb Z}_{p})$ we find
%\[{\mathbb M}(G_{i})=e_{i}H^{1}_{DR}({\mathcal A}/{\mathbb Z}_{p})=W({\mathbb F}_{p^2}).\]
%Therefore the height of $G_{i}$ is $2$ and {\bf Proposition 3.1} (in particular its proof) implies the claim.
%\begin{flushright}
%$\Box$
%\end{flushright}

\begin{prop} Let $J$ be an abelian variety of dimension $g$ over ${\mathbb F}_{p}$ so that
\[J[p^{\infty}]=G_{1}\times \cdots \times G_{g},\]
over ${\mathbb F}_{p}$, where $G_i$ is a $p$-divisible group defined over ${\mathbb F}_p$ with $\overline{G}_{i}\simeq G_{1,1}$ ($\forall i$). Then it is a product of  supersingular elliptic curve $\{E_{i}\}_{1\leq i \leq g}$ 
\[J=E_{1}\times  \cdots \times E_{g},\]
over ${\mathbb F}_{p}$. If moreover $p$ is greater than $3$,  $\Phi_{J,p}(t)=(t^2+p)^g$.
\end{prop}
{\bf Proof.} By (4) the $a$-number of $J$ is $g$ and, by {\bf Fact 2.4}, 
\[J\otimes_{{\mathbb F}_{p}}\overline{\mathbb F}_{p} \simeq E_{1}\times \cdots \times E_{g},\]
where $E_i$ is a supersingular elliptic curve ($\forall i$).
%we see that $J\otimes_{{\mathbb F}_{p}}\overline{\mathbb F}_{p}$ is isomorphic to a product of supersingular elliptic curves, $E_{1}\times \cdots \times E_{g}$. 
A simple consideration shows that, if necessary changing indices, we may assume that $G_{i}=E_{i}[p^{\infty}]$. Let $\phi$ be the $p$-th power Frobenius. Since $J$ is defined over ${\mathbb F}_{p}$ we have the diagram:
\[\begin{CD}
 E_{i}@>\text{$\phi$}>>  E_{i}^{\phi}\\
 @V\text{$\nu_{i}$}VV  @V\text{$\nu_{i}^{\phi}$}VV\\
  J  @>\text{$\phi$}>>   J,
 \end{CD}\]
 where $\nu_{i}$ is the imbedding.  Take the $p$-divisible groups and
 \[\begin{CD}
 E_{i}[p^{\infty}]@>\text{$\phi$}>>  E_{i}^{\phi}[p^{\infty}]\\
 @V\text{$\nu_{i}[p^{\infty}]$}VV  @V\text{$\nu_{i}^{\phi}[p^{\infty}]$}VV\\
  J[p^{\infty}]  @>\text{$\phi$}>>   J[p^{\infty}].
 \end{CD}\]
Since by the assumption all $\{E_{i}[p^{\infty}](=G_{i})\}_{i}$ and $\{\nu_{i}[p^{\infty}]\}_{i}$ are defined over ${\mathbb F}_{p}$, $E_{i}[p^{\infty}]=E_{i}^{\phi}[p^{\infty}]$ and $\nu_{i}[p^{\infty}]=\nu_{i}^{\phi}[p^{\infty}]$ ($\forall i$). Thus 
$E_{i}=E_{i}^{\phi}$ and $\nu_{i}=\nu_{i}^{\phi}$, which shows that each component $\{E_{i}\}_{i}$ and the product are defined over ${\mathbb F}_{p}$. The last claim follows from the well-known fact that the characteristic polynomial of $p$-power Frobenius of a supersingular elliptic curve over ${\mathbb F}_{p}$ is $t^{2}+p$ if $p\geq 5$ (\cite{Stichtenoth}).
\begin{flushright}
$\Box$
\end{flushright}
{\bf Theorem 1.3} follows from (11), {\bf Lemma 4.1} and {\bf Proposition 4.1}. \\

Let $F$ be a field of characteristic $0$ and $A$ an abelian variety over $F$ of dimension $g$. Fix a base $\{\omega_{1},\cdots,\omega_{g}\}$ of $H^{0}(A,\Omega^{1})$ over $F$ and we consider the faithful representation,
\[\rho : {\rm End}_{\overline{F}}(\overline{A}) \to M_{g}(\overline{F}),\]
defined by
\[\alpha^{*}\omega=\omega\cdot \rho(\alpha), \quad \omega=(\omega_{1},\cdots,\omega_{g}).\]
This is compatible with the action of ${\rm Gal}(\overline{F}/F)$ and the faithfulness of $\rho$ yields,
\begin{equation}{\rm End}_{F}(A)=\{ \alpha \in {\rm End}_{\overline{F}}(\overline{A})\,:\, \rho(\alpha)\in M_{g}(F)\}.\end{equation}
\begin{remark}
$\rho$ may be identified with the representation on the cotangent space of $A$ at the origin.
\end{remark}
\noindent{\bf Proof of Theorem 1.4}. (1) and (2) are consequences of {\bf Theorem 1.2} and {\bf Theorem 1.1}, respectively. Let us show (3). We first claim that (b) implies (a). Take a base $\{\omega_{1},\cdots,\omega_{g}\}$ of $H^{0}({\rm Jac}(C), \Omega^{1})$ over ${\mathbb Q}$ and consider the representation $\rho$. The assumption that the action of ${\mathcal O}_{K}$ is defined over $K$ implies $\rho({\mathcal O}_{K}) \subset M_{g}(K)$ by (12). Since $\rho$ is compatible with action of Galois group, take the invariant part of the complex conjugation and $\rho({\mathcal O}_{K_{0}})\subset M_{g}(K_{0})$. 
%\[\rho : {\mathcal O}_{K_{0}} \to M_{g}(K_{0}).\]
Thus, by (12), ${\mathcal O}_{K_{0}}$ is contained in ${\rm End}_{K_{0}}({\rm Jac}(C)\otimes_{{\mathbb Q}}K_{0})$. Now we show that (a) implies the claim. Let ${\mathcal J}$ be the Neron's model of ${\rm Jac}(C)$ over ${\mathbb Z}_{p}$, which is an abelian scheme. Since $p$ completely splits in $K_{0}$, ${\mathbb Z}_{p} \simeq ({\mathcal O}_{K_{0}})_{{\frak P}_{i}}$ (here we have used the notation of {\bf Theorem 1.3}). Together with (12) this implies that the action of ${\mathcal O}_{K_{0}}$ on ${\mathcal J}$ is defined over ${\mathbb Z}_{p}$. In fact by the Neron's mapping property of ${\mathcal J}$ (\cite{B-L-R}, \$1 {\bf Proposition 8}) it is sufficient to show that the action of ${\mathcal O}_{K_{0}}$ on the generic fiber ${\rm Jac}(C)\otimes_{\mathbb Q}{\mathbb Q}_{p}$ is defined over ${\mathbb Q}_{p}$. 
% Since the action of ${\mathcal O}_{K_{0}}$ is defined over $K_{0}$, $\rho$ induces a homomorphism,
%\[\rho : {\mathcal O}_{K_{0}} \to M_{g}(K_{0}).\]
Since $(K_{0})_{{\frak P}_{i}}\simeq {\mathbb Q}_{p}$ the image $\rho({\mathcal O}_{K_{0}})$ is contained in $M_{g}((K_{0})_{{\frak P}_{i}})\simeq M_{g}({\mathbb Q}_{p})$ and (12) shows that ${\mathcal O}_{K_{0}}\subset {\rm End}_{{\mathbb Q}_{p}}({\rm Jac}(C)\otimes_{\mathbb Q}{\mathbb Q}_{p})$. Take the special fiber and ${\rm Jac}(C_{p})$ satisfies the assumption of {\bf Theorem 1.3}.  
%By the assumption the Jacobian variety ${\rm Jac}(C_{p})$ of the reduction $C_{p}$ is the special fiber of  is an abelian scheme ${\mathcal J}$ over ${\mathbb Z}_{p}$ with CM and, 
Therefore ${\rm Jac}(C_{p})$ is a product of supersingular elliptic curves defined over ${\mathbb F}_{p}$.  Suppose that $p$ is greater than $3$. Since $\Phi_{{\rm Jac}(C_{p}),p}(t)=(t^2+p)^g$ the eigenvalues of the action of the $p$-th power Frobenius on $H^{1}_{et}(\overline{C}_{p}, {\mathbb Q}_{l})$ are $\{\sqrt{-p}, -\sqrt{-p}\}$. Use the Grothendieck-Lefschetz trace formula (\cite{Milne} {\bf Theorem 12.3}) and 
\[|C_{p}({\mathbb F}_{p^2})|=1+p^2-{\rm Tr}[Fr_{p^2} : H^{1}_{et}(\overline{C}_{p}, {\mathbb Q}_{l})]=1+p^2+2gp.\]
\begin{flushright}
$\Box$
\end{flushright}
%\noindent{\bf Proof of Theorem 1.4}.  (1) and (2) are consequences of {\bf Corollary 3.1}. {\bf Theorem 1.3} and Hasse-Weil's formula yield (3).
%. The assumption implies that the completion of ${\mathcal O}_{K_{0}}$ by any ${\frak P}_{i}$ is isomorphic to ${\mathbb Z}_{p}$ and ${\mathcal O}_{K_{0}}$ is contained in ${\mathbb Z}_{p}$.  Thus the Neron's minimal model of the Jacobian over ${\mathbb Z}_{p}$ satisfies the assumption of {\bf Theorem 1.3}. Therefore we see that the characteristic polynomial of $p$-th power Frobenius of the reduction of $C$ at $p$ is $(t^{2}+p)^{g}$. The last claim follows from the Hasse-Weil's formula.
%
%\begin{flushright}
%$\Box$
%\end{flushright}

%\begin{thm}
%$A$ is a product of  super-singular elliptic curve $\{E_{i}\}_{1\leq i \leq g}$: 
%\[A=E_{1}\times  \cdots \times E_{g},\]
%defined over ${\mathbb F}_{p}$. In particular $\Phi_{A,p}(t)=(t^2+p)^g$.
%\end{thm}

\noindent{\bf Proof of Corollaries}. We first show {\bf Corollary 1.1}. $K_{0}$ is ${\mathbb Q}(\zeta_{N}+\zeta_{N}^{-1})$ and the sequence
\[1 \to {\rm Gal}(K/K_{0})\simeq\{\pm 1\} \to {\rm Gal}(K/{\mathbb Q})\simeq ({\mathbb Z}/(N))^{\times} \to  {\rm Gal}(K_{0}/{\mathbb Q})\to 1,\]
shows that $p\equiv -1 (N)$ iff $p$ completely splits in $K_{0}$ and every prime factor of $p$ in $K_{0}$ remains prime in $K$. On the other hand $p\equiv 1 (N)$ iff $p$ completely splits in $K$. An existence of a positive integer $h$ with $p^{h}\equiv -1(N)$ implies that every prime factor of $p$ in $K_{0}$ remains prime in $K$. Now the desired claims follow from {\bf Theorem 1.4}. In the case of {\bf Corollary 1.2} observe $K_{0}={\mathbb Q}(\zeta_{M}+\zeta_{M}^{-1})$ and 
\[{\rm Gal}(K/{\mathbb Q})\simeq {\rm Gal}({\mathbb Q}(\zeta_{M}+\zeta_{M}^{-1})/{\mathbb Q})\times {\rm Gal}({\mathbb Q}(\zeta_{d})/{\mathbb Q}).\]
This is isomorphic to $(({\mathbb Z}/(M))^{\times}/{\pm 1}) \times \{\pm 1\}$ and the proof is similar.
\begin{flushright}
$\Box$
\end{flushright}

\section{Examples}
%Let $F$ be a field of characteristic $0$ and $A$ an abelian variety over $F$ of dimension $g$. We fix a base $\{\omega_{1},\cdots,\omega_{g}\}$ of $H^{0}(A,\Omega^{1})$ over $F$. Consider a faithful linear map
%\[\rho : {\rm End}_{\overline{F}}(\overline{A}) \to M_{g}(\overline{F}),\]
%defined to be
%\[\alpha^{*}\omega=\omega\cdot \rho(\alpha), \quad \omega=(\omega_{1},\cdots,\omega_{g}).\]
%This is compatible with the action of ${\rm Gal}(\overline{F}/F)$ and
%\begin{equation}{\rm End}_{F}(A)=\{ \alpha \in {\rm End}_{\overline{F}}(\overline{A})\,:\, \rho(\alpha)\in M_{g}(F)\}.\end{equation}
%A proper smooth curve $C$ defined over a number field $F$ of genus $g$ is mentioned as {\it a CM-curve} if the endomorphism ring of the Jacobian ${\rm Jac}(C_{\overline{F}})$ contains the integer ring of a CM-field $K$ so that $[K:{\mathbb Q}]=2g$. We say a finite prime $v$ of $F$ is {\it good} if the reduction $C_{v}$ is smooth and mention it {\it ordinary} (resp. {\it super-singular}, {\it super-special}) if the Jacobian has a corresponding reduction.\\

\noindent{\bf Example 5.1.}(\cite{Lang})
Let $l$ be an odd prime and we define a curve $C(l)$ to be the smooth projective model of
\[y^{l}=x(1-x),\] 
over ${\mathbb Q}$. The genus is $(l-1)/2$ and an $l$-th primitive root of unity $\zeta_{l}$  acts by
\[\zeta_{l}(x)=x,\quad \zeta_{l}(y)=\zeta_{l} y.\]
Since it is defined over $K:={\mathbb Q}(\zeta_{l})$, so is the action of ${\mathbb Z}[\zeta_{l}]$ on ${\rm Jac}(C(l))$. 
%Let $\{\omega_{1}\cdots, \omega_{(l-1)/1}\}$ be a base of $H^{0}(C(l),\Omega^{1})$ over ${\mathbb Q}$ and consider the representation $\rho$ in the previous section. The action of $\zeta_{l}$ on ${\rm Jac}(C(l))$ is defined over $K={\mathbb Q}(\zeta_{l})$ and $\rho(\zeta_{l})\in M_{(l-1)/2}(K)$ (In fact one may take them as eigenvectors of the action of $\zeta_{l}$ so that $(\zeta_{l})^{*}(\omega_{i})=\zeta_{l}^{\mu_{i}}$, $0\leq \mu_{i}\leq l-1$, \cite{Lang} \$1 {\bf Theorem 7.1}). Since $\rho$ is compatible with the Galois action $\rho(\zeta_{l}+\zeta_{l}^{-1})$ is contained in $M_{(l-1)/2}(K_{0})$ where $K_{0}={\mathbb Q}(\zeta_{l}+\zeta_{l}^{-1})$, and by (12) the action of ${\mathcal O}_{K_{0}}(={\mathbb Z}[\zeta_{l}+\zeta_{l}^{-1}])$ on ${\rm Jac}(C(l))$ is defined over $K_{0}$. 
Thus $C(l)$ satisfies the assumption {\bf Theorem 1.4}. It is also known that there is a ${\mathbb Q}$-rational base $\{\omega_{1}\cdots, \omega_{(l-1)/1}\}$ of $H^{0}(C(l),\Omega^{1})$ satisfying 
\[(\zeta_{l})^{*}(\omega_{i})=\zeta_{l}^{\mu_{i}}, \quad 0\leq \mu_{i}\leq l-1,\]
(\cite{Lang} \$1 {\bf Theorem 7.1}).
%%%%%%%%%
Let $p$ be a good prime so that $p^{h}\equiv -1(l)$ for a certain positive integer $h$. Set $q=p^{h}$ and let $F_{q+1}$ denote the Fermat curve,
\[X^{q+1}+Y^{q+1}=1.\]
By
\[x=X^{q+1},\quad y=(XY)^{\frac{q+1}{l}},\]
we have a surjective morphism defined over ${\mathbb Q}$,
\[\pi : F_{q+1} \to C(l),\]
and $H^{1}_{et}(C(l)_{p}\otimes_{{\mathbb F}_{p}}{\overline{\mathbb F}_{p}}, {\mathbb Q}_{l})$ is a Galois submodule of $H^{1}_{et}
((F_{q+1})_{p}\otimes_{{\mathbb F}_{p}}{\overline{\mathbb F}_{p}}, {\mathbb Q}_{l})$. On the other hand $|F_{q+1}({\mathbb F}_{p^2})|$ attains the Hasse-Weil upper bound $1+p^2+2gp$ and 
\[\Phi_{F_{q+1},p}(t)=(t^{2}+p)^{g},\]
where $g=q(q-1)/2$ is the genus of $F_{q+1}$ (\cite{Stichtenoth}{\bf Example 6.3.6}). Thus $\Phi_{C(l),p}(t)=(t^{2}+p)^{(l-1)/2}$ and 
\[|C(l)_{p}({\mathbb F}_{p^2})|=1+p^{2}+(l-1)p,\]
if there is a positive integer $h$ such that $p^{h}\equiv -1(l)$. {\bf Corollary 1.1} recovers this observation if $h=1$ but otherwise it only states that $C(l)$ has a supersingular reduction. Therefore {\bf Corollary 1.1} is only a sufficient condition for a prime $p$ at which the reduction of a CM-curve is ${\mathbb F}_{p^2}$-maximal, i.e. the number of ${\mathbb F}_{p^2}$-points attains the Hasse-Weil upper bound.\\

\noindent{\bf Example 5.2.}(\cite{Geemen-Koike}{\bf Example 2.3})
Let us consider a curve
\[C : y^{3}=x(x^{7}+1),\]
that has automorphisms
\[(x,y) \stackrel{\zeta_{7}}\mapsto (\zeta_{7}^{3}x, \zeta_{7}y), \quad  (x,y) \stackrel{\zeta_{3}}\mapsto (x, \zeta_{3}y),\]
and the involution,
\[\tau : (x,y) \mapsto (\frac{1}{x}, \frac{y}{x^{3}}).\]
The quotient $C/<\tau>$ has a smooth model
\[X : 2y^{3}=x^{4}-2\cdot 7^{2}x^{2}+2^{3}\cdot 7^{2}x-7^{3},\]
of genus $3$ and the endomorphism ring of ${\rm Jac}(X)$ is the integer ring ${\mathcal O}_{K}$ of $K:={\mathbb Q}(\zeta_{7}+\zeta_{7}^{-1},\zeta_{3})$.  Therefore $X$ is a CM-curve. A base of $H^{0}(X,\Omega^{1})$ is given by
\[\frac{(1-x)dx}{y},\quad \frac{(1-x^{4})dx}{y^{2}}, \quad \frac{(x-x^{3})dx}{y^{2}},\]
which are eigenvectors of the action of $\zeta_{7}+\zeta_{7}^{-1}$ whose eigenvalues are contained in $K_{0}={\mathbb Q}(\zeta_{7}+\zeta_{7}^{-1})$. Hence by (12) the action of ${\mathcal O}_{K_0}$ on ${\rm Jac}(X)$ is defined over $K_{0}$ and $X$ satisfies the assumption of {\bf Theorem 1.4}. By {\bf Corollary 1.2} we see that a good prime $p$ is supersingular if $p\equiv -1(3)$. Moreover it is superspecial if $p\equiv \pm 1(7)$ and $p\equiv -1(3)$. On the other hand, using \cite{Geemen-Koike}, one can determine $\Phi_{X,p}(t)$. In fact they have shown that ${\rm Jac}(C)$ is isogeneous to ${\rm Jac}(X)^{2}\times E$ where $E$ is a CM-elliptic curve whose defining equation is $x^{3}=y(y+1)$.

%In fact let $C$ be a curve over ${\mathbb Q}$ defined by
%\[C : y^{3}=x(x^{7}+1).\]
%Then Geemen, Koike and Weng 
\begin{fact}(\cite{Geemen-Koike}{\bf Corollary 3.1} and {\bf Corollary 3.2})
\begin{enumerate}
\item For $p \equiv 2, 5, 11, 17 \,(21)$, 
\[\Phi_{C,p}(t)=(t^{4}-pt^{2}+p^{2})^{2}(t^{2}+p)^{3}.\]

%\item For $p \equiv 2,11 (21)$, 
%\[\Phi_{C,p}(t)=(t^{4}-pt^{2}+p^{2})^{2}(t^{2}+p)^{3}.\]

\item For $p \equiv 8,20 \,(21)$, 
\[\Phi_{C,p}(t)=(t^{2}+p)^{7}.\]

\end{enumerate}
\end{fact}
Since $X$ has genus $3$ the degree of $\Phi_{X,p}(t)$ should be six and therefore
\begin{itemize}
\item If $p \equiv 2, 5, 11, 17 \,(21)$, 
\[\Phi_{X,p}(t)=(t^{4}-pt^{2}+p^{2})(t^{2}+p),\]
\item If $p \equiv 8,20 \,(21)$, 
\[\Phi_{X,p}(t)=(t^{2}+p)^{3}.\]
\end{itemize}
These coincide with the above results of {\bf Corollary 1.2} (see the discussion after {\bf Theorem 1.1}). Moreover in the second case ${\rm Jac}(X_{p})$ is a product of supersingular elliptic curves over ${\mathbb F}_{p}$ and $|X_{p}({\mathbb F}_{p^2})|=1+p^2+6p$.\\

\noindent{\bf Example 5.3.}(\cite{Geemen-Koike}{\bf Example 2.1})
Remember that the  {\it $n$-th Chebyshev polynomial} $U_{n}$ is defined by the recursive relation,
\[U_{n+1}(x)=xU_{n}(x)-U_{n-1}(x), \quad U_{0}(x)=2,\quad U_{1}(x)=x.\]
For a prime $l \geq 5$ we define a curve $X_{l}$ as
\[X(l) : y^{2}=U_{l}(x).\]
The genus is $(l-1)/2$ and it is the quotient of a hyperelliptic curve
\[Y(l): y^{2}=x(x^{2l}+1),\]
by the involution
\[\tau : (x,y) \mapsto (\frac{1}{x}, \frac{y}{x^{l+1}}).\]
Set $K={\mathbb Q}(\zeta_{l}+\zeta_{l}^{-1},\zeta_{4})$ and $K_{0}={\mathbb Q}(\zeta_{l}+\zeta_{l}^{-1})$. The automorphisms of $Y(l)$:
\[(x,y) \stackrel{\zeta_{l}}\mapsto (\zeta_{l}^{2}x, \zeta_{l}y), \quad  (x,y) \stackrel{\zeta_{4}}\mapsto (-x, \zeta_{4}y),\]
induce an action of ${\mathcal O}_{K}$ on ${\rm Jac}(X(l))$. The authors have shown that $H^{0}(X(l),\Omega^{1})$ has a ${\mathbb Q}$-rational base $\{\omega_{1},\cdots,\omega_{(l-1)/2}\}$ satisfying
\[(\zeta_{l}+\zeta_{l}^{-1})^{*}\omega_{i}=(\zeta_{l}^{i}+\zeta_{l}^{-i})\omega_{i},\quad 1\leq i \leq \frac{l-1}{2}.\]
By (12) this implies that the action of ${\mathcal O}_{K_{0}}$ on the Jacobian is defined over $K_{0}$. Thus $X(l)$ satisfies the assumption of {\bf Theorem 1.4}. {\bf Corollary 1.2} shows that a good prime $p$ satisfying $p\equiv -1(4)$ is supersingular. If moreover $p\equiv -1(4)$ and $p\equiv \pm 1(l)$, ${\rm Jac}(X(l)_{p})$ is a product of supersingular elliptic curves over ${\mathbb F}_{p}$ and $|X(l)_{p}({\mathbb F}_{p^{2}})|$ attains the Hasse-Weil upper bound $1+p^2+(l-1)p$. \\

%%%%%%%%%%%%
%\noindent{\bf Acknowledgement}. This research is partially supported by JSPS grants Kiban(C)22540068.
%

%%%%%%%%%%%%%%

%Address : Ken-ichi Sugiyama, \\
%Department of Mathematics and Informatics, Faculty of Science, Chiba University, \\
%1-33 Yayoi-cho Inage-ku,Chiba 263-8522, Japan,\\
%e-mail address : sugiyama@math.s.chiba-u.ac.jp

\end{document}